\documentclass[11pt]{amsart} \usepackage{latexsym, amssymb, stmaryrd,amsthm, amsmath, amscd}
\usepackage[T1]{fontenc}

\DeclareTextSymbol{\thh}{T1}{254}

\newtheorem{thm}{Theorem}[section]
\newtheorem{lemma}[thm]{Lemma}
\newtheorem{prop}[thm]{Proposition}
\newtheorem{cor}[thm]{Corollary}

\theoremstyle{definition}
\newtheorem{df}[thm]{Definition}
\newtheorem{rmk}[thm]{Remark}

\newtheorem{question}[thm]{Question}

\newcommand{\Z}{\mathbb{Z}}

\newcommand{\curly}[1]{\mathcal{#1}}
\newcommand{\A}{\curly{A}}

\newcommand{\D}{\curly{D}}

\newcommand{\m}{\mathcal{M}}
\newcommand{\N}{\mathcal{N}}

\newcommand{\cN}{\mathcal{N}}
\newcommand{\cM}{\mathcal{M}}

\newcommand{\la}{\curly{L}}

\makeatletter

\def\indsym#1#2{%
  \setbox0=\hbox{$\m@th#1x$}%
  \kern\wd0%
  \hbox to 0pt{\hss$\m@th#1\mid$\hbox to 0pt{$\m@th#1^{#2}$}\hss}%
  \lower.9\ht0\hbox to 0pt{\hss$\m@th#1\smile$\hss}%
  \kern\wd0}

\def\nindsym#1#2{%
  \setbox0=\hbox{$\m@th#1x$}%
  \kern\wd0%
  \hbox to 0pt{\hss$\m@th#1\not$\kern1.4\wd0\hss}
  \hbox to 0pt{\hss$\m@th#1\mid$\hbox to 0pt{$\m@th#1^{\,#2}$}\hss}%
  \lower.9\ht0\hbox to 0pt{\hss$\m@th#1\smile$\hss}%
  \kern\wd0}

\def\dotminussym#1#2{%
  \setbox0=\hbox{$\m@th#1-$}%
  \kern.5\wd0%
  \hbox to 0pt{\hss\hbox{$\m@th#1-$}\hss}%
  \raise.6\ht0\hbox to 0pt{\hss$\m@th#1.$\hss}%
  \kern.5\wd0}
\newcommand{\dotminus}{\mathbin{\mathpalette\dotminussym{}}}

\def \r { {\mathbb R} }
\def \<{\langle}
\def \>{\rangle}
\def \n {\mathbb N}

\def \*Z {{{^*}\Z}}

\def \((  {(\!(}
\def \)) {)\!)}

\def \int{\operatorname{int}}
\numberwithin{equation}{section}

\def \Aut{\operatorname{Aut}}

\def \D{\mathfrak{D}}
\def \c{\mathbb{C}}

\def \la{\mathcal{L}}
\def \D{\operatorname{D}}
\def \Age{\operatorname{Age}}
\def \Sym{\operatorname{Sym}}
\def \d{\mathbb{D}}
\def \ms{\operatorname{ms}}

\allowdisplaybreaks[2]

\begin{document}

\title[Approximate Herbrand and Definable Functions]{An Approximate Herbrand's Theorem and Definable Functions in metric structures}

\author{Isaac Goldbring}
\thanks{The author's work was partially supported by NSF grant DMS-1007144.}
\address {University of California, Los Angeles, Department of Mathematics, 520 Portola Plaza, Box 951555, Los Angeles, CA 90095-1555, USA}
\email{isaac@math.ucla.edu}
\urladdr{www.math.ucla.edu/~isaac}

\begin{abstract}
We develop a version of Herbrand's theorem for continuous logic and use it to prove that definable functions in infinite-dimensional Hilbert spaces are piecewise approximable by affine functions.  We obtain similar results for definable functions in Hilbert spaces expanded by a group of generic unitary operators and Hilbert spaces expanded by a generic subspace.  We also show how Herbrand's theorem can be used to characterize definable functions in some absolutely ubiquitous structures from classical logic.
\end{abstract}
\maketitle

\section{Introduction}

The main motivation for this paper comes from the study of definable functions in metric structures; this study was initiated by the author in \cite{GoldUry}, where a study of the definable functions in Urysohn's metric space was undertaken, and continued in \cite{Gold}, where the definable linear operators in (infinite-dimensional) Hilbert spaces were characterized.  However, lacking any understanding of arbitrary definable functions in Hilbert spaces, we conjectured that they were, in some sense, ``piecewise affine'' in analogy with the classical case of an infinite vector space over a division ring.  In unpublished lecture notes by van den Dries on motivic integration \cite{vdD}, we came upon a proof of the piecewise affineness of definable functions in such vector spaces using the following classical theorem of Herbrand:

\begin{thm}[Herbrand \cite{herb}]\label{classicalherb}
Suppose that $\la$ is a first-order signature and $T$ is a universal $\la$-theory with quantifier elimination.  Let $\varphi(\vec x,\vec y)$ be a formula, where $\vec x=(x_1,\ldots,x_m)$, $y=(y_1,\ldots,y_n)$, $m\geq 1$.  Then there are $\la$-terms 
$$t_{11}(\vec x),\ldots,t_{1n}(\vec x),\ldots,t_{k1}(\vec x),\ldots,t_{kn}(\vec x), \quad (k\in \n^{>0})$$ such that 
$$T\models \forall \vec x\forall \vec y\left(\varphi(\vec x,\vec y)\to \bigvee_{i=1}^k \varphi(\vec x,t_{i1}(\vec x),\ldots,t_{in}(\vec x)\right).$$
\end{thm}

Although this theorem is not immediately applicable to the case of an infinite vector space $V$ over a division ring (for the axioms expressing that $V$ is infinite are existential), Herbrand's theorem does apply to the theory of $V$ with constants added for names of elements of $V$.  Since terms in this extended language name affine functions, we get the aforementioned characterization of definable functions in $V$.  (According to van den Dries, this use of Herbrand's theorem is well-known and often used.)  Although Theorem \ref{classicalherb} has an easy model-theoretic proof using compactness, we should remark that the result was first established using proof-theoretic techniques; see \cite{Buss1} and \cite{Buss2} for more on the history of Herbrand's result.

In this paper, we prove a version of Herbrand's theorem for continuous logic (Theorem \ref{contherb} and Corollary \ref{herbT} below) and use it to characterize definable functions in Hilbert spaces and some of their generic expansions, proving, in the case of pure Hilbert spaces, that definable functions are ``piecewise approximable by affine functions.''  Along the way, we note that this method works whenever $T$ is a $\exists \forall$-axiomatizable theory with quantifier elimination.  In particular, we show that one can use Herbrand's theorem to understand definable functions in some absolutely ubiquitous structures from classical logic.

We assume that the reader is familiar with the basic definitions of continuous logic; otherwise, they can consult the survey article \cite{BBHU}.  

The author would like to thank Vinicius C.L., Aleksander Ivanov, and Dugald Macpherson for helpful discussions concerning this work and Matthias Aschenbrenner for pointing out the paper \cite{Lac} on absolutely ubiquitous structures.

\section{Herbrand's Theorem in Continuous Logic}

\

\noindent In this section, we let $\la$ denote an arbitrary continuous signature.  We will use the following abuse of notation:  whenever $\Delta$ is a set of closed $\la$-conditions and $\sigma$ is an $\la$-sentence, we write $\sigma \in \Delta$ to indicate that the condition ``$\sigma=0$'' belongs to $\Delta$. 
\begin{df}
Suppose that $\Delta$ is a set of closed $\la$-conditions.  
\begin{enumerate}
\item We say that $\Delta$ is \emph{closed under min} if whenever $\sigma_1,\ldots,\sigma_n$ are sentences with $\sigma_i\in \Delta$ for each $i$, then $\min_{1\leq i\leq n}\sigma_i\in \Delta$.
\item We say that $\Delta$ is \emph{closed under weakening} if whenever $\sigma\in \Delta$, then $\sigma\dotminus r\in \Delta$ for every $r\in [0,1]$.
\end{enumerate}
\end{df}

The following lemma is in a similar spirit to Lemma 3.4 of \cite{Usvy}; the classical version, whose proof we mimic, can be found in \cite{CK}.

\begin{lemma}\label{axiom}
Suppose that $T$ is a satisfiable $\la$-theory and $\Delta$ is a set of closed $\la$-conditions that is closed under min and weakening.  Then the following are equivalent:
\begin{enumerate}
\item $T$ is axiomatizable by a collection of conditions $\Gamma\subseteq \Delta$;
\item For all $\la$-structures $\m$ and $\N$ satisfying $\m\models T$ and $\sigma^\N=0$ for all $\sigma\in \Delta$ with $\sigma^\m=0$, we have $\N\models T$.
\end{enumerate}
\end{lemma}

\begin{proof}
Clearly $(1)\Rightarrow (2)$, so we need to prove $(2)\Rightarrow(1)$.  Consider the set $\Gamma=\{``\sigma=0\text{''} \ : \ \sigma\in \Delta \text{ and }T\models \sigma=0\}$.  We claim that $\Gamma$ axiomatizes $T$.  Suppose $\N\models \Gamma$.  Let $$\Sigma=\{``\delta \geq \frac{r}{2}\text{''} \ : \ \N\models \delta=r, \ r>0, \ \delta\in \Delta\}.$$  We claim that $T\cup \Sigma$ is consistent.  Suppose otherwise.  Then there are $\delta_1,\ldots,\delta_k$, $r_1,\ldots,r_k$ such that $T\models \min_{1\leq i\leq k} (\delta_i\dotminus \frac{r_i}{2})=0$.  Since $\Delta$ is closed under min and weakening, we have that $\min_{1\leq i\leq k} (\delta_i\dotminus \frac{r_i}{2})\in \Gamma$, so $\N\models \min_{1\leq i\leq k} (\delta_i\dotminus \frac{r_i}{2})=0$, which is a contradiction to the fact that $\delta_i^\N=r_i$ for each $i$.  Let $\m\models T\cup\Sigma$.  Now suppose that $\sigma\in \Delta$ and $\sigma^\m=0$.  Then $\sigma^\N=0$, else $``\sigma\geq \frac{r}{2}\text{''}\in \Sigma$ for some $r>0$, contradicting $\sigma^\m=0$.  By (2), we have $\N\models T$.
\end{proof}

Given an $\la$-structure $\m$, let $\D(\m)$ be the set of closed $\la(\m)$-conditions of the form $\sigma=0$, where $\sigma$ is an atomic $\la(\m)$ sentence and $\sigma^\m=0$; this is just the \emph{atomic diagram of $\m$}.  The following lemma is proved just as in classical logic.

\begin{lemma}
If $\N\models D(\m)$, then the $\la$-reduct of $\N$ contains a substructure isomorphic to $\m$.
\end{lemma}

Let us call a sentence $\sigma$ \emph{universal} if it is of the form $\sup_{\vec x}\varphi(\vec x)$, where $\varphi$ is quantifier-free.  Let us call a closed condition $``\sigma=0$'' \emph{universal} if $\sigma$ is universal.  We call a closed condition ``$\sigma=0$'' \emph{almost universal} if there is a universal sentence $\tau$ such that, in every $\la$-structure $\m$, we have $\sigma^\m=0$ if and only if $\tau^\m=0$.
   
\begin{lemma}
The set of almost universal conditions is closed under min and weakening.
\end{lemma}

\begin{proof}
Suppose that $\sigma=0$ and $\tau=0$ are almost universal conditions.  Suppose that $\sigma=0$ is equivalent to $\sup_{\vec x}\sigma'(\vec x)=0$ and $\tau=0$ is equivalent to $\sup_{\vec y}\tau'(\vec y)=0$, with $\sigma'$, $\tau'$ quantifier-free and $\vec x$, $\vec y$ disjoint tuples of distinct variables.  Then $\min(\sigma,\tau)=0$ is equivalent to $\sup_{\vec x}\sup_{\vec y}(\min(\sigma'(\vec x),\tau'(\vec y)))=0$.  Similarly, the condition $\sigma\dotminus r=0$ is equivalent to $\sup_{\vec x}(\sigma'(\vec x)\dotminus r)=0$.
\end{proof}

\noindent If $\Gamma$ is a set of closed $\la$-conditions, we set $$\Gamma^+:=\{``\sigma\leq \frac{1}{n}\text{''} \ : \ \sigma\in \Gamma, n\geq 1\}.$$

We say that $T$ has a \emph{universal axiomatization} if $T$ is axiomatizable by a set of universal conditions.  Clearly if $T$ is axiomatizable by a set of almost universal conditions, then $T$ has a universal axiomatization.

\begin{cor}\label{univ}
The following are equivalent:
\begin{enumerate}
\item $T$ has a universal axiomatization;
\item For any $\m\models T$ and substructure $\N$ of $\m$, we have $\N\models T$.
\end{enumerate}
\end{cor}

\begin{proof}
Clearly (1) implies (2), so we prove that (2) implies (1).  We use the criterion developed in Lemma \ref{axiom} applied to the set of almost universal conditions.  Suppose that $\m\models T$ and for all almost universal conditions $``\sigma=0$'', we have $\sigma^\m=0$ implies $\sigma^\N=0$.  We want $\N\models T$.  Let $T'=T\cup \D(\N)^+$.  We claim that $T'$ is satisfiable.  Fix atomic $\la(\N)$-sentences $\sigma_1(\vec b),\ldots,\sigma_n(\vec b)$ such that $\sigma_i^\N(\vec b)=0$.  Then $\N\models \inf_{\vec x}\max(\sigma_i(\vec x))=0$.  Suppose, towards a contradiction, that $\m\not\models \inf_{\vec x}\max(\sigma_i(\vec x))=0.$  Then there is $r\in (0,1]$ such that $\m\models \sup_{\vec x}(r\dotminus \max(\sigma_i(\vec x)))=0$.  By assumption, we have $\N\models \sup_{\vec x}(r\dotminus \max(\sigma_i(\vec x)))=0$, which is a contradiction.  Consequently, for any $k\geq 1$, there is $\vec a\in M$ such that $\m\models \max(\sigma_i(\vec a))\leq \frac{1}{k}$.  It follows by compactness that $T'$ is satisfiable.  Let $\A'\models T'$ and let $\A$ be the $\la$-reduct of $\A'$.  Then $\A\models T$ and $\N$ is (isomorphic to) a substructure of $\A$, whence $\N\models T$.
\end{proof}

\begin{df}
Suppose that $\m$ is an $\la$-structure and $A\subseteq M$.  Let $\langle A\rangle_0$ be the $\la$-prestructure generated by $A$.  Then the closure of $\langle A\rangle_0$ in $M$ is the completion of $\langle A\rangle_0$, whence a substructure of $\m$, called the \emph{substructure of $\m$ generated by $A$}.
\end{df}

\noindent By Theorem 3.5 of \cite{BBHU}, any $\la$-formula $\varphi(\vec x)$ has a modulus of uniform continuity $\Delta_\varphi:(0,1]\to (0,1]$, that is, for any $\la$-structure $\m$, any $\epsilon>0$, and any tuples $\vec a,\vec b$ from $M$, if $d(\vec a,\vec b)<\Delta_\varphi(\epsilon)$, then $|\varphi^\m(\vec a)-\varphi^\m(\vec b)|\leq \epsilon$.

\begin{thm}[Continuous Herbrand Theorem]\label{contherb}
Suppose that $T$ is a complete $\la$-theory with quantifier elimination that admits a universal axiomatization.  Let $\vec x=(x_1,\ldots,x_m)$ and $\vec y=(y_1,\ldots,y_n)$.  Then for any formula $\varphi(\vec x,\vec y)$ and any $\epsilon>0$, there are $\la$-terms $$t_{11}(\vec x),\ldots,t_{1n}(\vec x),\ldots,t_{k1}(\vec x),\ldots,t_{kn}(\vec x) \quad (k\in \n^{>0})$$ such that, for any $\m\models T$ and any $\vec a\in M^m$, if $\m\models \inf_{\vec y} \varphi(\vec a,\vec y)=0$, then $$\m \models \min_{1\leq i \leq k}\varphi(\vec a,t_{i1}(\vec a),\ldots,t_{in}(\vec a))\leq \epsilon.$$
\end{thm}

\begin{proof}
Consider the set of closed $\la$-conditions $\Gamma(\vec x)$ given by $$\{\inf_{\vec y}\varphi(\vec x,\vec y)=0\}\cup \{\varphi(\vec x,t_1(\vec x),\ldots,t_n(\vec x))\geq 2\epsilon \ : \ t_1(\vec x),\ldots,t_n(\vec x)\  \la\text{-terms}\}.$$  By compactness, it is enough to prove that $\Gamma$ is unsatisfiable.  Suppose, towards a contradiction, that $\m\models \Gamma(\vec a)$, where $\vec a=(a_1,\ldots,a_m)\in M^m$.  Fix $\delta\in (0,1]$ such that $\delta<\frac{\epsilon}{3}$.  Let $\chi(\vec x)$ be a quantifier-free $\la$-formula such that $T\models \sup_{\vec x}\left (|\inf_{\vec y}\varphi(\vec x,\vec y)-\chi(\vec x)|\dotminus \delta\right)=0$.  Then $\chi^\m(\vec a)\leq \delta$.  Let $\N$ be the substructure of $\m$ generated by $\{a_1,\ldots,a_m\}$.  Then since $\chi(\vec x)$ is quantifier-free, we have $\chi^\N(\vec a)\leq \delta$.  Since $\N\models T$, we have $\N\models \inf_{\vec y}\varphi(\vec a,\vec y)\leq 2\delta$.  Thus, there is $\vec c\in N^n$ such that $\varphi^\N(\vec a,\vec c)\leq 3\delta$.  Now let $t_i(\vec x)$ be a term so that $d(t_i(\vec a),\vec c_i)<\Delta_\varphi(\delta)$, whence $\varphi^\N(\vec a,t_1(\vec a),\ldots,t_n(\vec a))\leq 4\delta$.  Let $\theta(\vec x,\vec y)$ be a quantifier-free $\la$-formula so that $T\models \sup_{\vec x, \vec y}\left(|\varphi(\vec x,\vec y)-\theta(\vec x,\vec y)|\dotminus \delta\right)=0$.  Then $\theta^\N(\vec a,t_1(\vec a),\ldots,t_n(\vec a))\leq 5\delta$, whence $\theta^\m(\vec a,t_1(\vec a),\ldots,t_n(\vec a))\leq 5\delta$ and hence $\varphi^\m(\vec a,t_1(\vec a),\ldots,t_n(\vec a))\leq 6\delta$.  Since $6\delta<2\epsilon$, this is a contradiction to the fact that $\m\models \Gamma(\vec a)$.
\end{proof}

The following rephrasing of the previous theorem more closely resembles the usual statement of Herbrand's theorem.

\begin{cor}\label{herbT}
Suppose that $T$ is a complete $\la$-theory with quantifier elimination that admits a universal axiomatization.  Let $\vec x=(x_1,\ldots,x_m)$ and $\vec y=(y_1,\ldots,y_n)$.  Then for any formula $\varphi(\vec x,\vec y)$ and any $\epsilon>0$, there are $\la$-terms $$t_{11}(\vec x),\ldots,t_{1n}(\vec x),\ldots,t_{k1}(\vec x),\ldots,t_{kn}(\vec x) \quad (k\in \n^{>0})$$ and an increasing continuous function $\alpha:[0,1]\to [0,1]$ satisfying $\alpha(0)=0$ such that $$T\models \sup_{\vec x}(( \min_{1\leq i \leq k}\varphi(\vec x,t_{i1}(\vec x),\ldots,t_{in}(\vec x))\dotminus \epsilon)\dotminus \alpha(\inf_{\vec y}(\varphi(\vec x,\vec y)))=0.$$
\end{cor}

\begin{proof}
This is immediate from the preceding theorem and Proposition 7.15 of \cite{BBHU}.
\end{proof}
%We now let $H$ be an infinite-dimensional Hilbert space and $\la$ be the 1-sorted language suitable for studying unit balls of Hilbert spaces.  Let $T$ be the $\la(H)$-theory of $H$.  First observe that if $K\models T$, then $K\models IHS$ and $H\subseteq K$.  Since IHS admits QE, we have $H\preceq K$.  Observe that $T$ has a universal axiomatization.  Indeed, HS+$\D(H)$ is a universal axiomatization of $T$.  Also, $T$ admits QE because IHS admits QE.  Consequently, we can apply Herbrand's theorem to $T$.

\section{Primitive theories with QE}

In this short section, $\la$ continues to denote an arbitrary (continuous) signature and $T$ denotes an $\la$-theory.
\begin{df}
Following \cite{Lac} (in the classical setting), we say that $T$ is \emph{primitive} if there exists sets of closed $\la$-conditions $\Gamma$ and $\Delta$, where $\Gamma$ consists of universal conditions and $\Delta$ consists of existential conditions, such that $\Gamma\cup \Delta$ axiomatizes $T$.
\end{df}

\begin{rmk}
In classical logic, it is mentioned in \cite{Lac} that $T$ is primitive if and only if:  whenever $\cM_0,\cM_1\models T$ and $\cM_0\subseteq \N\subseteq \cM_1$, then $\N\models T$.  It is also mentioned in \cite{Lac} that $T$ is $\exists \forall$-axiomatizable if and only if:  whenever $\cM_0,\cM_1\models T$, $\cM_0\preceq \cM_1$, and $\cM_0\subseteq \N\subseteq \cM_1$, then $\N\models T$.  It follows that for model-complete theories $T$, $T$ is primitive if and only if $T$ is $\exists \forall$-axiomatizable.  An interesting example of a model-complete $\exists \forall$-theory is Example 3 of \cite{Lac2}.
\end{rmk}

\begin{prop}\label{prim}
Suppose that $T$ is a complete, model-complete primitive $\la$-theory.  Let $\cM\models T$ and let $T_\cM$ be the $\la(\m)$-theory of $\cM$.  Then $T_\cM$ is universally axiomatizable.  Moreover, $T_\cM$ has quantifier elimination if $T$ does.
\end{prop}

\begin{proof}
Let $\Gamma$ be a set of universal sentences and $\Delta$ a set of existential sentences such that $\Gamma\cup \Delta$ axiomatizes $T$.  In order to prove that $T_\cM$ has a universal axiomatization, it suffices to prove that $T_\cM$ is axiomatized by $\Gamma\cup \operatorname{D}(\cM)$.  Suppose that $\cN\models \Gamma\cup \operatorname{D}(\cM)$.  Then $\cM$ is a substructure of $\cN$.  Now any axiom from $\Delta$ is true in $\cN$ since it is witnessed by things in $\cM$.  Consequently, $\cN\models T$, whence $\cN\models T_\cM$ by model-completeness of $T$.  The moreover statement is clear.
%Now let $\sigma(\vec m)$ be an $\la(M)$-sentence so that $\cM\models \sigma(\vec m)=0$.  We would like $\cN\models \sigma(\vec m)=0$.  Let $\epsilon>0$.  Let $\psi(\vec x)$ be a quantifier-free $\la$-formula such that $T\models \sup_x|\sigma(\vec x)-\psi(\vec x)|\leq \epsilon$.  We have that $\cM\models \psi(\vec m)\leq \epsilon$, whence $\cN\models \psi(\vec m)\leq \epsilon$; since $\cN\models T$, we have $\cN\models \sigma(\vec m)\leq 2\epsilon$.  Letting $\epsilon\to 0$, we get that $\cN\models \sigma(\vec m)=0$.
\end{proof}

We will meet some examples of (classical and continuous) primitive theories with quantifier elimination in the next section.

The following proposition explains how we use Herbrand's theorem in connection with definable functions.
\begin{prop}\label{deffunc}
Suppose that $T$ is primitive and admits quantifier elimination.  Suppose $\cM\models T$ and $f:M^n\to M$ is a definable function.  Then for any $\epsilon>0$, there are $\la(M)$-terms $t_1(\vec x),\ldots,t_k(\vec x)$ such that: for all $\vec a\in M^n$, there is $i\in \{1,\ldots,k\}$ with $d(f(\vec a),t_i(\vec a))\leq \epsilon$.
\end{prop}
\begin{proof}
Fix $\epsilon>0$. Let $\varphi(\vec x,y)$ be an $\la(M)$-formula such that $$|d(f(\vec a),b)-\varphi^\cM(\vec a,b)|\leq \frac{\epsilon}{3}$$ for all $\vec a\in M^n$ and $b\in M$.  By Herbrand's theorem applied to $T_\m$ (which is applicable by Proposition \ref{prim}), there are $\la(M)$-terms $t_1(\vec x),\ldots,t_k(\vec x)$ such that, for all $\vec a\in M^n$, if $\cM\models \inf_y(\varphi(\vec a,y)\dotminus\frac{\epsilon}{3})=0$, then $$\cM\models (\varphi(\vec a,t_i(\vec a))\dotminus \frac{\epsilon}{3})\leq \frac{\epsilon}{3}$$ for some $i\in \{1,\ldots,k\}$.  Notice that the antecedent of the preceding conditional statement holds since $\varphi^\m(\vec a,f(\vec a))\leq \frac{\epsilon}{3}$.  Consequently, for every $\vec a\in M^n$, there is $i\in \{1,\ldots,k\}$ such that $d(f(\vec a),t_i(\vec a))\leq \epsilon$. 
\end{proof}

\begin{rmk}
Fix a definable function $f:M^n\to M$.  Fix $\epsilon>0$ and let the $\la(\m)$-terms $t_1(\vec x),\ldots,t_k(\vec x)$ be as in the conclusion of the previous proposition.  Suppose that $\cM\preceq \cN$ and $f:N^n\to N$ is the natural extension of $f$ to a definable function in $\cN$.  Then, for every $\vec a\in N^n$, there is $i\in \{1,\ldots,k\}$ such that $d(f(\vec a),t_i(\vec a))\leq \epsilon$.  Indeed, repeat the proof of the preceding proposition, using Corollary \ref{herbT} instead of Theorem \ref{contherb}.
\end{rmk}

\section{Applications}

In this section, we present some (classical and continuous) primitive theories with quantifier-elimination and use Proposition \ref{deffunc} above to understand the definable functions in models of these theories.

\subsection{Infinite-dimensional Hilbert spaces and some of their generic expansions}

\

\

\noindent In this subsection, we suppose that $\mathbb{K}\in \{\r,\c\}$ and we set $$\d:=\{\lambda \in \mathbb{K} \ : \ |\lambda|\leq 1\}.$$  Also, $\la$ denotes the (1-sorted) continuous signature for unit balls of $\mathbb{K}$-Hilbert spaces.  More specifically, $\la$ contains:
\begin{itemize}
\item a constant symbol $0$;
\item a binary function symbol $f_{\alpha,\beta}$ for every $\alpha,\beta\in \d$ with $|\alpha|+|\beta|\leq 1$;
\item a binary predicate symbol $\langle \cdot,\cdot \rangle$ that takes values in $[-1,1]$.
 \end{itemize}

If $H$ is a $\mathbb{K}$-Hilbert space, the unit ball of $H$, $B_1(H)$, is naturally an $\la$-structure, where $0$ is interpreted as the zero vector of $H$, $f_{\alpha,\beta}$ is interpreted as the function $(x,y)\mapsto \alpha x+\beta y$, and $\langle \cdot ,\cdot\rangle$ is interpreted as the inner product of $H$.  For sake of readability, we often write $H$ instead of $B_1(H)$ when speaking of this way of treating $B_1(H)$ as an $\la$-structure.  

Let $T$ be the $\la$-theory of (the unit ball of) an infinite-dimensional $\mathbb{K}$-Hilbert space.  Then $T$ is primitive as the Hilbert space axioms are universal and the axioms for infinite-dimensionality are existential.  We must remark that we cannot work in the many-sorted setting for Hilbert spaces (as in \cite{Gold}) because the axioms for the inclusion mappings are $\forall \exists$; indeed, for $n\leq m$, one must declare that the inclusion mapping $I_{n,m}:B_n(H)\to B_m(H)$ is onto the set of elements of $B_m(H)$ of norm at most $n$.

In the rest of this subsection, $H\models T$ and $H^*$ is an elementary extension of $H$.  In order to make any sense of Proposition \ref{deffunc} in this context, we must first understand $\la(H)$-terms.

\begin{lemma}\label{term}
If $t(x)$ is an $\la(H)$-term, then there are $\lambda\in \d$ and $v\in B_1(H)$ so that $t(a)=\lambda a+v$ for all $a\in B_1(H)$.
\end{lemma}

\begin{proof}
One proves this by induction on the complexity of $t(x)$, the base case being immediate.  Now suppose that $t_i(x)=\lambda_ix+v_i$ for $i=1,2$ and $\alpha,\beta$ are so that $|\alpha|+|\beta|\leq 1$.  Then $$f_{\alpha,\beta}(t_1(a),t_2(a))=\alpha t_1(a)+\beta t_2(a)=(\alpha \lambda_1+\beta \lambda_2)a+(\alpha v_1+\beta v_2).$$  It remains to observe that $|\alpha \lambda_1+\beta\lambda_2|\leq 1$.
\end{proof}

\begin{cor}\label{affine}
Let $f:H\to H$ be definable.  Then given $\epsilon>0$, there are $\lambda_1,\ldots,\lambda_k\in \mathbb{D}$ and $v_1,\ldots,v_k\in B_1(H)$ such that, for all $a\in B_1(H^*)$, there is $i\in \{1,\ldots,k\}$ with $d(f(a),\lambda_ia+v_i)\leq \epsilon$.
\end{cor}

Fix $a\in B_1(H^*)$.  Then there are sequences $(\lambda_n)$ from $ \mathbb{D}$ and $(v_n)$ from $B_1(H)$ with $\lambda_na+v_n\to f(a)$ as $n\to \infty$.  By taking subsequences, we may suppose that $\lambda_n\to \lambda\in \d$.  It then follows that $(v_n)$ is a Cauchy sequence in $B_1(H)$, whence $v_n\to v\in B_1(H)$.  It follows that $f(a)=\lambda a+v$.  We have just proven the following result:

\begin{cor}
For any $a\in B_1(H^*)$, there are $\lambda\in \mathbb{D}$ and $v\in B_1(H)$ such that $f(a)=\lambda a+v$.
\end{cor}

\begin{cor}
Suppose that $H^*$ is $\omega_1$-saturated and $f(H^\perp)\subseteq H^\perp$. Fix $\epsilon>0$ and let $\lambda_1,\ldots,\lambda_m$ be a finite $\epsilon$-net for $\d$.   Then there is a finite-dimensional subspace $K$ of $H$ such that, for all $a\in B_1(H^*)\cap K^\perp$, there is $i\in\{1,\ldots,m\}$ such that $d(f(a),\lambda_ia)<\epsilon$.    
\end{cor}

\begin{proof}
Let $a\in B_1(H^*)\cap H^\perp$.  Take $\lambda\in \d$ and $v\in B_1(H)$ such that $f(a)=\lambda a+v$.  Then $$0=\langle f(a),v\rangle =\langle \lambda a+v,v\rangle=\langle v,v\rangle.$$  Thus, $f(a)=\lambda a$.  Let $(a_n)$ be an orthonormal basis for $H$.  Then the following set of conditions is unsatisfiable in $H^*$:
$$\{\langle x,a_n\rangle=0 \ : \ n<\omega\} \cup \{d(f(x),\lambda_i x)\geq \epsilon \ : \ i=1,\ldots,m\}.$$  By saturation, there is $n<\omega$ such that, setting $K:=\operatorname{span}(a_1,\ldots,a_n)$, we have $d(f(x),\lambda_ix)<\epsilon$ for all $x\in B_1(H^*)\cap K^\perp$. 
\end{proof}
%\begin{cor}
%There are $(\lambda_i \ : \ i<2^{\aleph_0})$ and $(v_i \ : \ i<2^{\aleph_0})$ from $B_1(H)$ such that, for any elementary extension $B_1(H^*)$ of $B_1(H)$, and for any $a\in B_1(H^*)$, there is $i<2^{\aleph_0}$ such that $f(a)=\lambda_ia+v_i$.
%\end{cor}

%Of course the preceding corollary has nontrivial content only when $B_1(H)^*$ is not separable.

How does Corollary \ref{affine} relate to functions definable in the many-sorted language for Hilbert spaces considered in \cite{Gold}?  In order to elucidate this, we first clarify how the syntax of continuous logic works in the case that the predicates take values in intervals other than $[0,1]$.  (This is omitted in the survey \cite{BBHU} and was communicated to me by Ward Henson.)  Let $\la'$ be a many-sorted (continuous) signature with sort set $S$.  In particular, one associates to each predicate symbol $P$ of $\la$ a closed, bounded interval $I_P$ in $\r$.  Then one also associates to each formula $\varphi$ a closed, bounded interval $I_\varphi$ in $\r$ as follows:
\begin{itemize}
\item Given two terms $t_1(\vec x)$ and $t_2(\vec x)$ of arity $(s_1,\ldots,s_n,s_{n+1})$, the formula $\varphi(\vec x)=d(t_1(\vec x),t_2(\vec x))$ is an atomic formula with $I_\varphi:=[0,N]$, where $N$ is the bound on the metric of sort $s_{n+1}$.
\item If $P$ is a predicate symbol of arity $(s_1,\ldots,s_n)$ and $t_1(\vec x),\ldots,t_n(\vec x)$ are terms such that $t_i$ takes values in sort $s_i$, then the formula $\varphi(\vec x)=P(t_1(\vec x),\ldots,t_n(\vec x))$ is an atomic formula with $I_\varphi:=I_P$.
\item Suppose that $\varphi_1(\vec x),\ldots,\varphi_n(\vec x)$ are formulae with associated intervals $I_{\varphi_1},\ldots,I_{\varphi_n}$.  Suppose that $u$ is a continuous function with domain $I_{\varphi_1}\times \cdots \times I_{\varphi_n}$ and range $I$, a closed, bounded interval in $\r$.  Then $\varphi(\vec x)=u(\varphi_1(\vec x),\ldots,\varphi_n(\vec x))$ is a formula with $I_\varphi:=I$.
\item If $\varphi$ is a formula with associated interval $I_\varphi$, then $\psi=\sup_x \varphi$ is a formula with $I_\psi:=I_\varphi$.  Similarly for $\inf_x\varphi$.
\end{itemize}

For an interval $I=[a,b]\subseteq \r$ with $a<b$, define $u_I:I\to [0,1]$ by $u_I(x):=\frac{1}{b-a}(x-a)$.  Note that $u_I$ is a homeomorphism with inverse $u_I^{-1}(x)=a+(b-a)x$.

We let $\la_{\ms}$ denotes the many-sorted theory of Hilbert spaces used in \cite{Gold}.
\begin{lemma}
For any quantifier-free $\la_{\ms}$-formula $\varphi(\vec x)$, where $\vec x$ is a tuple of variables of sort $B_1(H)$, there is a quantifier-free $\la$-formula $\psi(\vec x)$ with $I_{\psi}=[0,1]$  such that $$H\models \sup_{\vec x}|u_{I_\varphi}(\varphi(\vec x))-\psi(\vec x)|=0.$$  In particular, when $I_\varphi=[0,1]$, we have $H\models \sup_x |\varphi(\vec x)-\psi(\vec x)|=0$.
\end{lemma}

\begin{proof}
The proof goes by induction on the complexity of $\varphi$, the main work taking place in the case when $\varphi$ is atomic, which involves a painful case distinction.  Let us illustrate the idea by considering  terms $t_i(x,y)=\lambda_ix+\mu_iy$ ($i=1,2$) where $|\lambda_i|,|\mu_i|\leq n$.  (In the general situation, terms can be much more complicated due to the number of variables and the inclusion maps.)  First suppose that $\varphi(x,y)=d(t_1(x,y),t_2(x,y))$.  Since each $t_i$ takes values in $B_{2n}$, we have $I_{\varphi}=[0,4n]$.  Then $I_{\varphi}(\varphi(x,y))=\frac{1}{4n}d(t_1(x,y),t_2(x,y))$.  Let $\psi(x,y)=\|\frac{\lambda_1-\lambda_2}{4n}x+\frac{\mu_1-\mu_2}{4n}y\|$.  Since $|\frac{\lambda_1-\lambda_2}{4n}|+|\frac{\mu_1-\mu_2}{4n}|\leq 1$, we have that $\psi$ is an $\la$-formula with $I_{\psi}=[0,1]$.  Clearly $\psi$ is as desired.  

Now  suppose that $\varphi(x,y)=\langle t_1(x,y),t_2(x,y)\rangle$.  Now $I_{\varphi}=[-4n^2,4n^2]$, so $u_{I_{\varphi}}(\varphi(x,y))=\frac{1}{8n^2}(\langle t_1(x,y),t_2(x,y)\rangle +4n^2)$.  This time, let $$\psi(x,y)=\frac{1}{2}\langle \frac{\lambda_1}{2n}x+\frac{\mu_1}{2n}y,\frac{\lambda_2}{2n}+\frac{\mu_2}{2n}y\rangle+\frac{1}{2}.$$  It is easily verified that this $\psi$ is as desired.

For the induction step, suppose that $\varphi=u(\varphi_1,\ldots,\varphi_n)$, where $$u:I_{\varphi_1}\times \cdots \times I_{\varphi_n}\to I_{\varphi}$$ is a surjective continuous function.  By the induction hypothesis, there are $\la$-formulae $\psi_i(x)$ ($i=1,\ldots,n$) with each $I_{\psi_i}=[0,1]$  such that $H\models \sup_{\vec x}|u_{I_{\varphi_i}}(\varphi_i(\vec x))-\psi_i(\vec x)|=0$.  Consider the $\la$-formula $$\psi(x)=u_{I_\varphi}(u(u_{I_{\varphi_1}}^{-1}(\psi_1(\vec x)),\ldots,u_{I_{\varphi_n}}^{-1}(\psi_n(\vec x)))).$$  It is clear that $H\models \sup_x|u_{\varphi}(\varphi(\vec x))-\psi(\vec x))|=0$.
\end{proof}

\begin{cor}
If $P:B_1(H)^n\to [0,1]$ is a uniformly continuous function, then $P$ is an $\la$-definable predicate if and only if $P$ is an $\la_{\operatorname{ms}}$-definable predicate
\end{cor}

\begin{proof}
This follows from the preceding corollary and the fact that the $\la_{ms}$-theory of $H$ admits quantifier-elimination.
\end{proof}

\begin{cor}
Suppose that $f:H\to H$ is an $\la_{\operatorname{ms}}$-definable function such that $f(B_1(H))\subseteq B_1(H)$.  Then $f|B_1(H)$ is an $\la$-definable function.
\end{cor}
\noindent The definition of an $\la_{\operatorname{ms}}$-definable function is given in \cite{Gold}.
%\begin{proof}
%Now suppose that $f:H\to H$ is a definable function in the many-sorted sense as defined in \cite{}.  By multiplying $f$ by a scalar if necessary, we may assume that $f(B_1(H))\subseteq B_1(H)$.  Also, by scaling the parameters if necessary, we may suppose that $f$ is $A$-definable, where $A\subseteq B_1(H)$.  We then claim that $f|B_1(H)$ is definable in the language $\mathcal{L}_n$.  The key point is that the many-sorted theory has QE, and the quantifier-free formulae in either language with variables from $B_n$ are the same.  Consequently, we have that the preceding theorems hold for our original class of definable functions.
%\end{proof}

\begin{rmk}
It follows from the preceding corollary and Corollary \ref{affine} that for any $\la_{\operatorname{ms}}$-definable function $f:H\to H$,  any $n\geq 1$, and any $\epsilon >0$, there are scalars $\lambda_1,\ldots,\lambda_k$ and vectors $v_1,\ldots,v_k \in B_{m(n,f)}(H)$ such that, for all $x\in B_n(H)$, there is $i\in \{1,\ldots,k\}$ with $d(f(x),\lambda_i x+v_i)\leq \epsilon$.  Using the main result of \cite{Gold}, we can give a different proof of this fact in the case that $f$ is linear.  Indeed, write $f=\lambda I+K$, where $K$ is a compact operator.  Let $\{v_1,\ldots,v_k\}$ be a finite $\epsilon$-net for $K(B_n(H))$.  Then for $a\in B_1(H)$, we have $d(K(a),v_i)\leq \epsilon$ for some $i\in \{1,\ldots,k\}$, whence $d(f(a),\lambda a+v_i)\leq \epsilon$.  (Notice here that $\lambda_i=\lambda$ for all $i$.)
\end{rmk}

\
%What we have said holds not just for $B_1(H)$, but rather for any $B_n(H)$.  

We now suppose that $\mathbb{K}=\c$ and set $\mathbb{S}^1:=\{\lambda\in \c \ : \ |\lambda|=1\}$.  We let $\la_U:=\la\cup \{U,U^{-1}\}$, where $U$ and $U^{-1}$ are both unary function symbols.  We let $T_U^\forall$ denote the $\la$-theory obtained from $T$ by adding (universal) axioms saying that $U$ is linear, preserves the inner product, and $U$ and $U^{-1}$ are inverses.  ($T_U$ axiomatizes the theory of an infinite-dimensional Hilbert space equipped with a unitary operator; one adds a symbol for $U^{-1}$ so as to avoid the $\forall \exists$ axiom stating that $U$ is onto.)  We add to $T^\forall_U$ the following axioms:
$$\inf_x[|\langle x,x\rangle -1|\dotplus d(Ux,\sigma x)|]=0,$$ where $\sigma$ ranges over a countable dense subset of $\mathbb{S}^1$.  (These axioms assert that the spectrum of $U$ is $\mathbb{S}^1$.)  Then $T_U$ is complete and admits quantifier elimination (see \cite{BUZ}); $T_U$ is the theory of infinite-dimensional Hilbert spaces equipped with a generic automorphism.  Since $T_U$ is primitive, we can once again apply Proposition \ref{deffunc}.

\begin{lemma}
If $t(x)$ is an $\la_U(H)$-term, then there are $l,m\in\Z$, $l\leq m$, $\alpha_l,\ldots,\alpha_m\in \mathbb{D}$ and a vector $v\in B_1(H)$ such that, for all $a\in B_1(H)$, we have $$t(a)=v+\sum_{j=l}^m\alpha_jU^j(a).$$
\end{lemma}

\begin{proof}
This is proved by induction on the complexity of $t(x)$ exactly as in Lemma \ref{term}.
\end{proof}

Suppose that $(H^*,U^*)$ is an elementary extension of $(H,U)$.

\begin{cor}
Suppose that $f:H\to H$ is an $\la_U$-definable function and $\epsilon>0$.  Then there are $l,m\in \Z$, $l\leq m$, $\lambda^1_l,\ldots,\lambda^1_m,\ldots,\lambda^k_l,\ldots,\lambda^k_m\in \mathbb{D}$, and $v_1,\ldots,v_k\in B_1(H)$, such that, for all $a\in B_1(H^*)$, there is $i\in \{1,\ldots,k\}$ such that $$d(f(x),v_i+\sum_{j=l}^m \alpha^i_jU^j(x))<\epsilon.$$
\end{cor}

%Using the compactness of $\d^\Z$ and arguing as in the proof of Corollary \ref{affine}, we obtain the following:

%\begin{cor}
%For any $a\in B_1(H^*)$, there are $(\lambda_k)\in \d^\Z$ and $v\in B_1(H)$ such that $\sum_{k\in \Z} \lambda_kU^k(a)$ exists and $f(a)=v+\sum_{k\in \Z}\lambda_kU^k(a)$.
%\end{cor}

One can generalize this situation as follows:  Let $G$ be a countable (discrete group) and let $\la_G$ be the language for Hilbert spaces as above augmented by unary function symbols $\tau_g$ for $g\in G$.  Let $T_G$ be the universal $\la_G$-theory of a unitary representation of $G$ on an infinite-dimensional Hilbert space.  (As above, the axiom $\sup_xd((\tau_g(\tau_{g^{-1}}(x)),x)=0$ allows us to assert that $\tau_g$ is onto without using a $\forall\exists$ axiom.) Let $\pi:G\to U(H)$ be a unitary representation of $G$ on an (infinite-dimensional) Hilbert space $H$ such that $(H,\pi)$ is an existentially closed model of $T_G$ (such an existentially closed model exists because $T_G$ is an inductive theory).  Let $\Sigma$ be the set of existential consequences of $(H,\pi)$.  Then it is shown in \cite{Ber} that $T_{GA}:=T_G\cup \Sigma$ axiomatizes the class of existentially closed models of $T_G$, whence is the model companion of $T_G$.  Moreover, since $T_G$ has the amalgamation property (see \cite{Ber}), it follows that $T_{GA}$ admits quantifier elimination.  As above, one can show that any $\la_G$ term $t(x)$ has the form $v+\sum_{i=1}^n\lambda_i g_ix$ for some $v\in B_1(H)$, some $\lambda_1,\ldots,\lambda_n\in \d$, and some $g_1,\ldots,g_n\in G$.  (Here we abuse notation and write $gx$ instead of $\tau_g(x)$.)  Consequently, we have:

\begin{cor}
Let $(H,\pi)$ be any model of $T_{GA}$ and let $f:H\to H$ be an $\la_G$-definable function.  Then, for any $\epsilon>0$, there are $v_1,\ldots,v_k\in B_1(H)$, scalars $\lambda^1_1,\ldots,\lambda^1_m,\ldots,\lambda^k_1,\ldots, \lambda^k_m\in \d$, and group elements $g_1,\ldots,g_k\in G$ such that, for all $a\in B_1(H^*)$, there is $i\in \{1,\ldots,k\}$ such that $$d(f(a),v_i+\sum_{j=1}^m\lambda^i_jg_j a)<\epsilon.$$  %Consequently, for any elementary extension $(H^*,\pi^*)$ of $(H,\pi)$ and any $a\in B_1(H^*)$, there are $(\lambda_n)\in \d^\n$ and $v\in B_1(H)$ such that $\sum_{n=1}^\infty \lambda_ng_na$ exists and $f(a)=v+\sum_{n=1}^\infty \lambda_ng_na$. 
\end{cor}

There is yet another expansion of Hilbert spaces that fits into this context. Let $\la_P:=\la\cup \{P\}$, where $P$ is a new unary predicate symbol.  We consider the theory $T_P$ obtained from the theory of infinite-dimensional Hilbert spaces obtained by adding the following axioms (the latter two are axiom schemes, including one such axiom for every $n\geq 1$):
\begin{itemize}
\item $P$ is linear;
\item $\sup_x d(P^2(x),P(x))=0$;
\item $\sup_{x,y}|\langle P(x),y\rangle-\langle x,P(y)\rangle|=0$;
\item $\inf_{v_1}\cdots \inf_{v_n} \max (\max_{i,j}|\langle v_i,v_j\rangle\dotminus \delta_{ij}|,\max_i d(P(v_i),v_i)))=0$;
\item $\inf_{v_1}\cdots \inf_{v_n} \max (\max_{i,j}|\langle v_i,v_j\rangle\dotminus \delta_{ij}|,\max_i d(P(v_i),0)))=0$.
\end{itemize}

The first three axioms say that $P$ is a projection operator on $H$ and the latter two axiom schemes say that $P(H)$ and $P(H)^\perp$ are infinite-dimensional.  Then $T_P$ is a complete theory with quantifier elimination (\cite{BV}); in fact, it is the theory of beautiful pairs of Hilbert spaces and its unique separable model is the Fraisse limit of the family of finite-dimensional Hilbert spaces equipped with projection operators.

Since $T_P$ is a primitive theory with quantifier elimination, we may use Proposition \ref{deffunc}.  Let $(H,P)$ be a model of $T_P$.  Then in $(H,P)$, all $\la$-terms $t(x)$ are easily seen to equivalent to terms be of the form $\alpha x+\beta P(x)+v$, where $\alpha,\beta\in \d$ and $v\in B_1(H)$.  Thus:

\begin{prop}
Let $f:B_1(H)\to B_1(H)$ be an $\la_P$-definable function.  Then for any $\epsilon>0$, there are $v_1,\ldots,v_k\in B_1(H)$ and $\alpha_1,\ldots,\alpha_k,\beta_1,\ldots,\beta_k\in \d$ such that, for all $a\in B_1(H)$, there is $i\in\{1,\ldots,k\}$ such that $$d(f(a),\alpha_ia+\beta_iP(a)+v_i)<\epsilon.$$  Consequently, for any elementary extension $(H^*,P^*)$ of $(H,P)$ and any $a\in B_1(H^*)$, there are $\alpha,\beta\in \d$ and $v\in B_1(H)$ such that $f(a)=\alpha a+\beta P^*(a)+v$.
\end{prop}

\subsection{Absolutely ubiquitous structures}

\

\

\noindent A source of primitive theories in classical logic comes from the notion of an \emph{absolutely ubiquitous} structure.  Suppose that $\la$ is a finite first-order signature and $\cM$ is a countable $\la$-structure.  Recall that $\cM$ is said to be \emph{locally finite} if every finitely generated substructure of $\cM$ is finite and $\cM$ is said to be \emph{uniformly locally finite} if there is a function $g:\n^{>0}\to\n^{>0}$ such that, for all $A\subseteq M$, if $|A|\leq n$, then $|\langle A\rangle|\leq g(n)$, where $\langle A\rangle$ denotes the substructure of $\cM$ generated by $A$.  Also recall that the \emph{age of $\cM$}, denoted $\Age(\cM)$, is the set of isomorphism classes of finitely generated substructures of $\cM$.  Finally, we say that $\cM$ is absolutely ubiquitous if:
\begin{enumerate}
\item $\cM$ is uniformly locally finite, and 
\item whenever $\cN$ is a countable, locally finite $\la$-structure with $\Age(\cM)=\Age(\N)$, then $\cM\cong \N$.  
\end{enumerate}
It follows immediately from the definition that if $\cM$ is an absolutely ubiquitous $\la$-structure and $T:=\operatorname{Th}(\cM)$, then $T$ is primitive and $\aleph_0$-categorical, whence model-complete (see also Lemma 2.1 of \cite{M2}).  Consequently, if $T$ has quantifier elimination, then $T$ meets the hypothesis of Proposition \ref{deffunc}.  

It is interesting to ask when an absolutely ubiquitous structure has quantifier elimination?  Note that an absolutely ubiquitous structure admits quantifier elimination if and only if it is ultrahomogeneous.  Thus, we can use the classifications of absolutely ubiquitous graphs \cite{M} and ultrahomogeneous countable graphs \cite{LW} to see that there are only two situations when a countable ultrahomogeneous graph is absolutely ubiquitous:
\begin{itemize}
\item a disjoint union of finitely many copies of the complete graph on $\aleph_0$ many vertices;
\item a $k$-partite graph, where each part is of size $\aleph_0$.
\end{itemize}

It follows from Proposition \ref{deffunc} that if $G$ is such a graph and $f:G^n\to G$ is a definable function, then there are vertices $g_1,\ldots,g_k\in G$ so that, for any $\vec a\in G^n$, we have $f(a)=a_i$ for some $i$ or $f(a)=g_j$ for some $j$.

It is interesting to note that in the case of absolutely ubiquitous structures in finite \emph{relational} signatures, we can always expand the language to ensure that we have quantifier elimination while maintaining absolute ubiquity.  To see this, suppose that $\cM$ is an $\la$-structure, where $\la$ is a finite relational (classical) signature.  We say that $\cM$ is \emph{finitely partitioned} if there is a finite partition $X_1,\ldots,X_n$ of $M$ such that $\Sym(X_1)\times\cdots \times \Sym(X_n)$ is a subgroup of $\Aut(\cM)$.  The main result of \cite{HM} states that $\cM$ is absolutely ubiquitous if and only if $\cM$ is finitely partitioned.  Suppose now that $\cM$ is absolutely ubiquitous.  Let $X_1,\ldots,X_n$ be a finite partition of $\cM$ witnessing that $\cM$ is finitely partitioned.  Consider the signature $\la':=\la\cup\{R_1,\ldots,R_n\}$, where $R_1,\ldots,R_n$ are new unary function symbols, and consider the expansion $\cM':=(\cM;X_1,\ldots,X_n)$ of $\cM$ to an $\la'$-structure.  Then $X_1,\ldots,X_n$ witness that $\cM'$ is finitely partitioned, whence $\cM'$ is absolutely ubiquitous.  However, we now have:

\begin{lemma}
$\cM'$ is ultrahomogeneous, whence $\operatorname{Th}(\cM')$ admits quantifier elimination.
\end{lemma}

\begin{proof}
Suppose that $A,B\subseteq M$ are finite and $f:A\to B$ is a partial automorphism of $\cM'$.  Then for any $i\in \{1,\ldots,n\}$, $f(A\cap X_i)\subseteq X_i$.  Extend $f$ to $\tilde{f}:M\to M$ so that $\tilde{f}|X_i\in \operatorname{Sym}(X_i)$ for each $i\in\{1,\ldots,n\}$.  Then by assumption, $\tilde{f}\in \Aut(\cM')$.
\end{proof}

\begin{cor}
Given any definable (in $\cM'$) function $f:M^n\to M$, there are elements $b_1,\ldots,b_m\in M$ so that, for all $\vec a\in M^n$, we have either $f(\vec a)=a_i$ for some $i\in \{1,\ldots,n\}$ or $f(\vec a)=b_j$ for some $j\in \{1,\ldots,m\}$.
\end{cor}

What about when the language has function symbols?  Here is an example from \cite{V}:  Let $\cM:=(\n^n,E_1,\ldots,E_n,f)$, where $E_i$ is the binary relation on $\n^n$ given by $E_i(\vec a,\vec b)\Leftrightarrow a_i=b_i$, and $f$ is the $n$-ary function on $\n^n$ given by $f(\vec a_1,\ldots, \vec a_n)=(a_{11},\ldots,a_{nn})$.  It is argued in \cite{V} that $\cM$ is an absolutely ubiquitous structure with quantifier elimination.

It is shown in \cite{M2} that if $G$ is an absolutely ubiquitous group (considered as a structure in the pure group language), then $G$ has a characteristic subgroup $A$ of finite index such that $A$ is a finite direct product of elementary abelian groups of infinite rank.  Conversely, if $
G$ is a countable group with a characteristic subgroup $A$ of index $m<\infty$ which is a finite direct product of elementary abelian groups of infinite rank such that either $G=A\times F$ for some group $F$ of cardinality $m$ or $m$ is relatively prime to the orders of elements of $A$, then $G$ is absolutely ubiquitous.  If the absolutely ubiquitous group $G$ admits quantifier elimination, then given any definable function $f:G^n\to G$, there is a tuple $\vec b$ from $G$ and words $w_1(\vec x,\vec b),\ldots,w_k(\vec x,\vec b)$, such that, for all $\vec a\in G^n$, there is $i\in \{1,\ldots,k\}$ such that $f(\vec a)=w_i(\vec a,\vec b)$.  

The question remains: which absolutely ubiquitous groups admit quantifier elimination?  It is easy to see that if $G$ itself is a finite direct product of elementary abelian groups of infinite rank, then $G$ is ultrahomogeneous, so admits quantifier elimination.  More generally:

\begin{prop} 
If $G=A\times F$, where $A$ is a finite direct product of elementary abelian groups of infinite rank, $F$ is a finite ultrahomogeneous group, and $\operatorname{gcd}(|a|,|b|)=1$ for all $a\in A$ and $b\in F$, then $G$ is ultrahomogeneous.
\end{prop}

\begin{proof}
Suppose that $\phi:B\to C$ is an isomorphism, where $B$ and $C$ are finite subgroups of $G$.  Let $A_1,F_1$ denote the projections of $B$ onto $A$ and $F$ respectively; note that $A_1$ and $F_1$ are finite subgroups of $A$ and $F$ respectively.  Next note that, for each $a\in A_1$, we have that $(a,1)\in B$.  Indeed, if $(a,b)\in B$, then choosing $n\in \n$ such that $|b|$ divides $n$ and $n\equiv 1\mod |a|$, we see that $(a,1)=(a,b)^n\in B$.  Likewise, for every $b\in F_1$, we have $(1,b)\in B$.  Now observe that, for all $(a,1)\in B$, there is $a'\in A$ such that $\phi(a,1)=(a',1)$.  Indeed, writing $\phi(a,1)=(a',b)$, we have $(1,1)=\phi(a,1)^{|a|}=(1,b^{|a|})$, whence $b=1$.  Similarly, for every $b\in F$, there is $b'\in F$ such that $\phi(1,b)=(1,b')$.  We can thus define $\phi':A_1\to A$ by $\phi'(a)=a'$, where $\phi(a,1)=(a',1)$; note that $\phi'$ is an isomorphism between finite subgroups of $A$, so can be lifted to an automorphism $\tilde{\phi}':A\to A$.  Likewise, one obtains a partial automorphism $\phi'':F_1\to F$ that can be lifted to an automorphism $\tilde{\phi}'':F\to F$.  Finally, $\tilde{\phi}:G\to G$ defined by $\tilde{\phi}(a,b)=(\tilde{\phi}'(a),\tilde{\phi''}(b))$ is an automorphism of $G$ extending $\phi$ 
\end{proof}

\begin{rmk}
The ultrahomogeneous finite groups are characterized in \cite{CF}.
\end{rmk}

\begin{question}
Given an absolutely ubiquitous group $G$, is there an extension $\la'$ of the language of groups by relation symbols and an expansion $\mathcal{G}$ of $G$ to an $\la'$-structure so that $\mathcal{G}$ admits quantifier elimination and is still absolutely ubiquitous (or at least has a primitive theory)?  If the answer to this question is positive, then definable functions in absolutely ubiquitous groups are piecewise given by words as mentioned above.
\end{question}

\end{document}